\documentclass[a4paper]{article}

\usepackage{amssymb}

\def\GAP{\textsf{GAP}}
\def\ATLAS{\textsc{Atlas}}
\def\Irr{{\rm Irr}}
\def\FF{{\mathbb{F}}}
\def\tthdump#1{#1}
\tthdump{}

\begin{document}

\tthdump{\title{{\GAP} computations needed in the proof of \\
\cite[Theorem~6.1~(ii)]{DNT}}}

\author{\textsc{Thomas Breuer} \\[0.3cm]
\textit{Lehrstuhl D f{\"u}r Mathematik} \\
\textit{RWTH, 52056 Aachen, Germany} \\[0.6cm]
\textsc{Klaus Lux} \\[0.3cm]
\textit{Department of Mathematics} \\
\textit{University of Arizona, Tucson, AZ 85721, USA}}

\date{September 19th, 2011}

\maketitle

\abstract{
This is a collection of example computations that are cited in the Appendix
of~\cite{DNT}.
In each case, the aim is to show that the extension of a given finite simple
group by an elementary abelian group of given rank has the property that
not all complex irreducible characters of the same degree are
Galois conjugate.

The purpose of this writeup is twofold.
On the one hand, the details of the computations are documented this way.
On the other hand, the {\GAP} code shown for the examples can be used as
test input for automatic checking of the data and the functions used.}

\textwidth16cm
\oddsidemargin0pt

\parskip 1ex plus 0.5ex minus 0.5ex
\parindent0pt

\tableofcontents


For the computations, we need some Brauer character tables from~\cite{JLPW95},
some generating matrices from~\cite{AGR},
and some functions from the {\GAP} system~\cite{GAP} and its packages
\verb|AtlasRep|~\cite{AtlasRep},
\verb|cohomolo|~\cite{cohomolo},
\verb|CTblLib|~\cite{CTblLib1.2}, and
\verb|TomLib|~\cite{TomLib}.

First we load the necessary {\GAP} packages.

\begin{verbatim}
    gap> LoadPackage( "AtlasRep", "1.5" );
    true
    gap> LoadPackage( "cohomolo", "1.6" );
    true
    gap> LoadPackage( "CTblLib", "1.2" );
    true
    gap> LoadPackage( "TomLib", "1.2.1" );
    true
\end{verbatim}

\section{$G/N \cong Sz(8)$ and $|N| = 2^{12}$}

The group $S = Sz(8)$ has exactly one irreducible $12$-dimensional module
over the field with two elements, up to isomorphism.
This module can be obtained from any of the three absolutely irreducible
$4$-dimensional $S$-modules in characteristic two,
by regarding it as a module over the prime field $\FF_2$.

\begin{verbatim}
    gap> p:= 2;;  d:= 12;;
    gap> t:= CharacterTable( "Sz(8)" ) mod p;
    BrauerTable( "Sz(8)", 2 )
    gap> irr:= Filtered( Irr( t ), x -> x[1] <= d );;
    gap> Display( t, rec( chars:= irr, powermap:= false, centralizers:= false ) );
    Sz(8)mod2
    
           1a 5a 7a 7b 7c 13a 13b 13c
    
    Y.1     1  1  1  1  1   1   1   1
    Y.2     4 -1  A  C  B   D   F   E
    Y.3     4 -1  B  A  C   E   D   F
    Y.4     4 -1  C  B  A   F   E   D
    
    A = E(7)^2+E(7)^3+E(7)^4+E(7)^5
    B = E(7)+E(7)^2+E(7)^5+E(7)^6
    C = E(7)+E(7)^3+E(7)^4+E(7)^6
    D = E(13)+E(13)^5+E(13)^8+E(13)^12
    E = E(13)^4+E(13)^6+E(13)^7+E(13)^9
    F = E(13)^2+E(13)^3+E(13)^10+E(13)^11
    gap> List( irr, x -> SizeOfFieldOfDefinition( x, p ) );
    [ 2, 8, 8, 8 ]
\end{verbatim}

First we construct the $12$-dimensional irreducible representation of $S$
over $\FF_2$,
using that the {\ATLAS} of Group Representations provides matrix generators
for $S$ in the $4$-dimensional representation over $\FF_8$.

\begin{verbatim}
    gap> info:= OneAtlasGeneratingSetInfo( "Sz(8)", Dimension, 4,
    >               Characteristic, p );
    rec( charactername := "4a", dim := 4, groupname := "Sz(8)", id := "a",
      identifier := [ "Sz(8)", [ "Sz8G1-f8r4aB0.m1", "Sz8G1-f8r4aB0.m2" ], 1, 8 ],
      repname := "Sz8G1-f8r4aB0", repnr := 17, ring := GF(2^3), size := 29120,
      standardization := 1, type := "matff" )
    gap> gens_dim4:= AtlasGenerators( info ).generators;;
    gap> b:= Basis( GF(8) );; 
    gap> gens_dim12:= List( gens_dim4, x -> BlownUpMatrix( b, x ) );;
\end{verbatim}

We claim that any extension of $S$ with the given module splits.

\begin{verbatim}
    gap> s:= AtlasGroup( "Sz(8)", IsPermGroup, true );;
    gap> chr:= CHR( s, p, 0, gens_dim12 );;
    gap> SecondCohomologyDimension( chr );
    0
\end{verbatim}

(The function \verb|CHR| takes as its arguments a permutation group,
the characteristic of the module, a finitely presented group (or zero),
and a list of matrices that define the module in the sense that they
correspond to the generators of the given permutation group.
Note that this condition is satisfied because the generators provided by
the {\ATLAS} of Group Representations are compatible.)

So it is enough to consider the semidirect product $G = 2^{12}\!:\!Sz(8)$.
%

The {\GAP} Character Table Library contains the ordinary character table
of $G$.
We check this as follows.
By the above cohomology result,
the group $G$ is uniquely determined, up to isomorphism, by the group order
and the property that $G$ has a minimal normal subgroup $N$
such that $G/N$ is a simple group isomorphic with $S$.
(Since $|G|/|S|$ is a power of two, $N$ is a $2$-group.
By the minimality condition, $N$ is elementary abelian and the action of $S$
on $N$ affords the desired $S$-module.
Note that the isomorphism type of a finite simple group is determined
by its character table.)

\begin{verbatim}
    gap> iso:= IsomorphismTypeInfoFiniteSimpleGroup( s );
    rec( name := "2B(2,8) = 2C(2,8) = Sz(8)", parameter := 8, series := "2B" )
    gap> names:= AllCharacterTableNames( Size, 2^12 * Size( s ) );;
    gap> cand:= List( names, CharacterTable );;
    gap> cand:= Filtered( cand,
    >      t -> ForAny( ClassPositionsOfMinimalNormalSubgroups( t ),
    >             n -> IsomorphismTypeInfoFiniteSimpleGroup( t / n ) = iso ) );
    [ CharacterTable( "2^12:Sz(8)" ) ]
\end{verbatim}

So we can easily check that $G$ has eight rational valued irreducibles
of the degree $455$ (or of the degree $3\,640$).

\begin{verbatim}
    gap> t:= cand[1];;
    gap> rationals:= Filtered( Irr( t ), x -> IsSubset( Integers, x ) );;
    gap> Collected( List( rationals, x -> x[1] ) );
    [ [ 1, 1 ], [ 64, 1 ], [ 91, 1 ], [ 455, 8 ], [ 3640, 8 ] ]
\end{verbatim}

\section{$G/N \cong M_{22}$ and $|N| = 2^{10}$}

The group $S = M_{22}$ has exactly two irreducible $10$-dimensional modules
over the field with two elements, up to isomorphism.
These modules are in fact absolutely irreducible.

\begin{verbatim}
    gap> p:= 2;;  d:= 10;;
    gap> t:= CharacterTable( "M22" ) mod p;
    BrauerTable( "M22", 2 )
    gap> irr:= Filtered( Irr( t ), x -> x[1] <= d );;
    gap> Display( t, rec( chars:= irr, powermap:= false, centralizers:= false ) );
    M22mod2
    
           1a 3a 5a 7a 7b 11a 11b
    
    Y.1     1  1  1  1  1   1   1
    Y.2    10  1  .  A /A  -1  -1
    Y.3    10  1  . /A  A  -1  -1
    
    A = E(7)+E(7)^2+E(7)^4
      = (-1+Sqrt(-7))/2 = b7
    gap> List( irr, x -> SizeOfFieldOfDefinition( x, p ) );
    [ 2, 2, 2 ]
\end{verbatim}

First we construct the two irreducible $10$-dimensional representations
of $S$ over $\FF_2$,
again using that the {\ATLAS} of Group Representations provides the matrix
generators in question.

\begin{verbatim}
    gap> info:= AllAtlasGeneratingSetInfos( "M22", Dimension, d,
    >               Characteristic, p );
    [ rec( charactername := "10a", dim := 10, groupname := "M22", id := "a",
          identifier := [ "M22", [ "M22G1-f2r10aB0.m1", "M22G1-f2r10aB0.m2" ], 1,
              2 ], repname := "M22G1-f2r10aB0", repnr := 13, ring := GF(2),
          size := 443520, standardization := 1, type := "matff" ),
      rec( charactername := "10b", dim := 10, groupname := "M22", id := "b",
          identifier := [ "M22", [ "M22G1-f2r10bB0.m1", "M22G1-f2r10bB0.m2" ], 1,
              2 ], repname := "M22G1-f2r10bB0", repnr := 14, ring := GF(2),
          size := 443520, standardization := 1, type := "matff" ) ]
    gap> gens:= List( info, r -> AtlasGenerators( r ).generators );;
\end{verbatim}

We claim that any extension of $S$ with any of the two given modules splits.

\begin{verbatim}
    gap> s:= AtlasGroup( "M22", IsPermGroup, true );;
    gap> chr:= CHR( s, p, 0, gens[1] );;
    gap> SecondCohomologyDimension( chr );
    0
    gap> chr:= CHR( s, p, 0, gens[2] );;
    gap> SecondCohomologyDimension( chr );
    0
\end{verbatim}

Again we see that it is enough to consider semidirect products
$G = 2^{10}\!:\!M_{22}$, but this time for the two nonisomorphic modules.
%

The {\GAP} Character Table Library contains the ordinary character tables
of the two groups in question.
We check this with the same approach as in the previous examples.

\begin{verbatim}
    gap> iso:= IsomorphismTypeInfoFiniteSimpleGroup( s );
    rec( name := "M(22)", series := "Spor" )
    gap> names:= AllCharacterTableNames( Size, 2^10 * Size( s ) );;
    gap> cand:= List( names, CharacterTable );;
    gap> cand:= Filtered( cand,
    >      t -> ForAny( ClassPositionsOfMinimalNormalSubgroups( t ),
    >             n -> IsomorphismTypeInfoFiniteSimpleGroup( t / n ) = iso ) );
    [ CharacterTable( "2^10:M22'" ), CharacterTable( "2^10:m22" ) ]
    gap> List( cand, NrConjugacyClasses );
    [ 47, 43 ]
\end{verbatim}

So we can easily check that in both cases,
$G$ has two rational valued irreducibles of the degree $1\,155$.

\begin{verbatim}
    gap> t:= cand[1];;
    gap> rationals:= Filtered( Irr( t ), x -> IsSubset( Integers, x ) );;
    gap> Collected( List( rationals, x -> x[1] ) );
    [ [ 1, 1 ], [ 21, 1 ], [ 22, 1 ], [ 55, 1 ], [ 99, 1 ], [ 154, 1 ], 
      [ 210, 1 ], [ 231, 3 ], [ 385, 1 ], [ 440, 1 ], [ 770, 5 ], [ 924, 2 ], 
      [ 1155, 2 ], [ 1386, 1 ], [ 1408, 1 ], [ 3080, 2 ], [ 3465, 4 ], 
      [ 4620, 2 ], [ 6930, 3 ], [ 9240, 1 ] ]
    gap> t:= cand[2];;
    gap> rationals:= Filtered( Irr( t ), x -> IsSubset( Integers, x ) );;
    gap> Collected( List( rationals, x -> x[1] ) );
    [ [ 1, 1 ], [ 21, 1 ], [ 55, 1 ], [ 77, 1 ], [ 99, 1 ], [ 154, 1 ], 
      [ 210, 1 ], [ 231, 1 ], [ 330, 1 ], [ 385, 3 ], [ 616, 2 ], [ 693, 1 ], 
      [ 770, 1 ], [ 1155, 2 ], [ 1980, 1 ], [ 2310, 4 ], [ 2640, 1 ], 
      [ 3465, 2 ], [ 4620, 1 ], [ 5544, 2 ], [ 6160, 1 ], [ 6930, 2 ], 
      [ 9856, 1 ] ]
\end{verbatim}

\section{$G/N \cong J_2$ and $|N| = 2^{12}$}

The group $S = J_2$ has exactly one irreducible $12$-dimensional module
over the field with two elements, up to isomorphism.
This module can be obtained from any of the two absolutely irreducible
$6$-dimensional $S$-modules in characteristic two,
by regarding it as a module over the prime field $\FF_2$.

\begin{verbatim}
    gap> p:= 2;;  d:= 12;;
    gap> t:= CharacterTable( "J2" ) mod p;
    BrauerTable( "J2", 2 )
    gap> irr:= Filtered( Irr( t ), x -> x[1] <= d );;
    gap> Display( t, rec( chars:= irr, powermap:= false, centralizers:= false ) );
    J2mod2
    
           1a 3a 3b 5a 5b 5c 5d 7a 15a 15b
    
    Y.1     1  1  1  1  1  1  1  1   1   1
    Y.2     6 -3  .  A *A  B *B -1   C  *C
    Y.3     6 -3  . *A  A *B  B -1  *C   C
    
    A = -2*E(5)-2*E(5)^4
      = 1-Sqrt(5) = 1-r5
    B = E(5)+2*E(5)^2+2*E(5)^3+E(5)^4
      = (-3-Sqrt(5))/2 = -2-b5
    C = E(5)+E(5)^4
      = (-1+Sqrt(5))/2 = b5
    gap> List( irr, x -> SizeOfFieldOfDefinition( x, p ) );
    [ 2, 4, 4 ]
\end{verbatim}

First we construct the irreducible $12$-dimensional representation of $S$
over $\FF_2$,
using that the {\ATLAS} of Group Representations provides matrix generators
for $S$ in the $6$-dimensional representation over $\FF_4$.

\begin{verbatim}
    gap> info:= OneAtlasGeneratingSetInfo( "J2", Dimension, 6,
    >               Characteristic, p );
    rec( charactername := "6a", dim := 6, groupname := "J2", id := "a",
      identifier := [ "J2", [ "J2G1-f4r6aB0.m1", "J2G1-f4r6aB0.m2" ], 1, 4 ],
      repname := "J2G1-f4r6aB0", repnr := 16, ring := GF(2^2), size := 604800,
      standardization := 1, type := "matff" )
    gap> gens_dim6:= AtlasGenerators( info ).generators;;
    gap> b:= Basis( GF(4) );;
    gap> gens_dim12:= List( gens_dim6, x -> BlownUpMatrix( b, x ) );;
\end{verbatim}

We claim that any extension of $S$ with the given module splits.

\begin{verbatim}
    gap> s:= AtlasGroup( "J2", IsPermGroup, true );;
    gap> chr:= CHR( s, p, 0, gens_dim12 );;
    gap> SecondCohomologyDimension( chr );
    0
\end{verbatim}

Again we see that it is enough to consider a semidirect product
$G = 2^{12}\!:\!J_2$.
%

The {\GAP} Character Table Library contains the ordinary character table
of $G$.
We check this with the same approach as in the previous examples.

\begin{verbatim}
    gap> iso:= IsomorphismTypeInfoFiniteSimpleGroup( s );
    rec( name := "HJ = J(2) = F(5-)", series := "Spor" )
    gap> names:= AllCharacterTableNames( Size, 2^12 * Size( s ) );;
    gap> cand:= List( names, CharacterTable );;
    gap> cand:= Filtered( cand,
    >      t -> ForAny( ClassPositionsOfMinimalNormalSubgroups( t ),
    >             n -> IsomorphismTypeInfoFiniteSimpleGroup( t / n ) = iso ) );
    [ CharacterTable( "2^12:J2" ) ]
\end{verbatim}

So we can easily check that $G$ has two rational valued irreducibles
of the degree $1\,575$.

\begin{verbatim}
    gap> t:= cand[1];;
    gap> rationals:= Filtered( Irr( t ), x -> IsSubset( Integers, x ) );;
    gap> Collected( List( rationals, x -> x[1] ) );
    [ [ 1, 1 ], [ 36, 1 ], [ 63, 1 ], [ 90, 1 ], [ 126, 1 ], [ 160, 1 ], 
      [ 175, 1 ], [ 225, 1 ], [ 288, 1 ], [ 300, 1 ], [ 336, 1 ], [ 1575, 2 ], 
      [ 2520, 4 ], [ 3150, 1 ], [ 4725, 6 ], [ 9450, 1 ], [ 10080, 4 ], 
      [ 12600, 4 ], [ 18900, 2 ] ]
\end{verbatim}

\section{$G/N \cong J_2$ and $|N| = 5^{14}$}

The group $S = J_2$ has exactly one irreducible $14$-dimensional module
over the field with $5$ elements, up to isomorphism.
This module is in fact absolutely irreducible.

\begin{verbatim}
    gap> p:= 5;;  d:= 14;;
    gap> t:= CharacterTable( "J2" ) mod p;
    BrauerTable( "J2", 5 )
    gap> irr:= Filtered( Irr( t ), x -> x[1] <= d );;
    gap> Display( t, rec( chars:= irr, powermap:= false, centralizers:= false ) );
    J2mod5
    
           1a 2a 2b 3a 3b 4a 6a 6b 7a 8a 12a
    
    Y.1     1  1  1  1  1  1  1  1  1  1   1
    Y.2    14 -2  2  5 -1  2  1 -1  .  .  -1
\end{verbatim}

In this case, we do not attempt to compute the complete character table of
$G$.
Instead, we show that $G/N$ has at least five regular orbits on the
dual space of $N$, and apply~\cite[Lemma~5.1~(i)]{DNT}.
(Note that $N$ is in fact self-dual.)

For that, we use {\GAP}'s table of marks of $S$.
The information stored for this table of marks allows us to compute,
for each class of subgroups $U$ of $S$, the numbers of orbits in the dual
space of $N$ for which contain the point stabilizers in $S$ are exactly
the conjugates of $U$.
The following {\GAP} function takes the table of marks \verb|tom| of $S$,
a list \verb|matgens| of matrices that describe the action of the generators of
$S$ on the vector space in question, and the size \verb|q| of its field of scalars.
The return value is a record with the components
\verb|fixed| (the vector of numbers of fixed points of the subgroups of $S$
on the dual of $N$),
\verb|decomp| (the numbers of orbits with the corresponding point stabilizers),
\verb|nonzeropos| (the positions of subgroups that occur as point stabilizers),
and \verb|staborders| (the list of orders of the subgroups that occur as
point stabilizers).

\begin{verbatim}
    gap> orbits_from_tom:= function( tom, matgens, q )
    >     local slp, fixed, idmat, i, rest, decomp, nonzeropos;
    > 
    >     slp:= StraightLineProgramsTom( tom );
    >     fixed:= [];
    >     idmat:= matgens[1]^0;
    >     for i in [ 1 .. Length( slp ) ] do
    >       if IsList( slp[i] ) then
    >         # Each subgroup generator has a program of its own.
    >         rest:= List( slp[i],
    >                      prg -> ResultOfStraightLineProgram( prg, gens ) );
    >       else
    >         # The subgroup generators are computed with one common program.
    >         rest:= ResultOfStraightLineProgram( slp[i], gens );
    >       fi;
    >       if IsEmpty( rest ) then
    >         # The subgroup is trivial.
    >         fixed[i]:= q^Length( idmat );
    >       else
    >         # Compute the intersection of fixed spaces of the transposed
    >         # matrices, since we act on Irr(N) not on N.
    >         fixed[i]:= q^Length( NullspaceMat( TransposedMat( Concatenation(
    >                        List( rest, x -> x - idmat ) ) ) ) );
    >       fi;
    >     od;
    > 
    >     decomp:= DecomposedFixedPointVector( tom, fixed );
    >     nonzeropos:= Filtered( [ 1 .. Length( decomp ) ],
    >                            i -> decomp[i] <> 0 );
    > 
    >     return rec( fixed:= fixed,
    >                 decomp:= decomp,
    >                 nonzeropos:= nonzeropos,
    >                 staborders:= OrdersTom( tom ){ nonzeropos },
    >               );
    > end;;
\end{verbatim}

Note that this function assumes that the generators of $S$ obtained from
the {\ATLAS} of Group Representations are compatible with the generators
from {\GAP}'s table of marks of $S$.
This fact can be read off from the \verb|true| value of the \verb|ATLAS| component
in the \verb|StandardGeneratorsInfo| value of the table of marks.

\begin{verbatim}
    gap> tom:= TableOfMarks( "J2" );
    TableOfMarks( "J2" )
    gap> StandardGeneratorsInfo( tom );
    [ rec( ATLAS := true, description := "|z|=10, z^5=a, |b|=3, |C(b)|=36, |ab|=7"
            , generators := "a, b", 
          script := [ [ 1, 10, 5 ], [ 2, 3 ], [ [ 2, 1 ], [ "|C(",, ")|" ], 36 ], 
              [ 1, 1, 2, 1, 7 ] ], standardization := 1 ) ]
\end{verbatim}

Alternatively, we can compute whether the generators are compatible,
as follows.

\begin{verbatim}
    gap> info:= OneAtlasGeneratingSetInfo( "J2", Dimension, d, Ring, GF(p) );
    rec( charactername := "14a", dim := 14, groupname := "J2", id := "", 
      identifier := [ "J2", [ "J2G1-f5r14B0.m1", "J2G1-f5r14B0.m2" ], 1, 5 ], 
      repname := "J2G1-f5r14B0", repnr := 19, ring := GF(5), size := 604800, 
      standardization := 1, type := "matff" )
    gap> gens:= AtlasGenerators( info ).generators;;
    gap> map:= GroupGeneralMappingByImages( UnderlyingGroup( tom ),
    >      Group( gens ), GeneratorsOfGroup( UnderlyingGroup( tom ) ), gens );;
    gap> IsGroupHomomorphism( map );
    true
\end{verbatim}

Now we are sure that we may apply the function \verb|orbits_from_tom|.

\begin{verbatim}
    gap> orbits_from_tom( tom, gens, p );
    rec( decomp := [ 8600, 0, 2512, 359, 10, 0, 0, 212, 5, 0, 0, 4, 0, 240, 16,
          10, 0, 0, 0, 0, 10, 0, 0, 0, 0, 2, 0, 0, 36, 0, 0, 0, 26, 0, 0, 0, 0,
          0, 0, 0, 20, 0, 10, 8, 0, 0, 0, 0, 0, 0, 0, 0, 0, 0, 0, 0, 0, 10, 0, 0,
          5, 0, 0, 0, 26, 0, 10, 0, 0, 0, 0, 10, 0, 0, 0, 0, 0, 0, 0, 0, 0, 0, 0,
          0, 0, 0, 0, 0, 10, 0, 0, 0, 2, 0, 0, 0, 0, 2, 4, 0, 0, 0, 0, 0, 4, 0,
          0, 0, 0, 0, 0, 0, 0, 0, 0, 16, 0, 0, 0, 0, 0, 0, 0, 0, 0, 2, 0, 0, 0,
          0, 0, 0, 0, 0, 0, 0, 0, 0, 4, 0, 0, 0, 4, 0, 0, 1 ],
      fixed := [ 6103515625, 15625, 390625, 390625, 625, 25, 3125, 3125, 625,
          625, 625, 625, 5, 3125, 125, 625, 25, 25, 125, 5, 125, 25, 125, 25, 25,
          25, 5, 125, 125, 125, 25, 25, 3125, 1, 1, 5, 5, 25, 5, 25, 125, 5, 25,
          25, 25, 25, 25, 25, 5, 25, 25, 5, 25, 5, 5, 5, 5, 25, 25, 1, 125, 1, 5,
          5, 125, 1, 25, 5, 25, 1, 5, 25, 5, 5, 25, 25, 5, 5, 5, 1, 5, 5, 1, 1,
          1, 5, 1, 25, 25, 25, 1, 5, 25, 5, 5, 1, 1, 125, 5, 5, 5, 25, 5, 5, 5,
          1, 1, 5, 5, 1, 5, 1, 5, 1, 1, 25, 5, 5, 1, 1, 1, 1, 5, 1, 1, 25, 1, 1,
          5, 1, 1, 5, 1, 5, 1, 1, 5, 1, 5, 1, 1, 1, 5, 1, 1, 1 ],
      nonzeropos := [ 1, 3, 4, 5, 8, 9, 12, 14, 15, 16, 21, 26, 29, 33, 41, 43,
          44, 58, 61, 65, 67, 72, 89, 93, 98, 99, 105, 116, 126, 139, 143, 146 ],
      staborders := [ 1, 2, 3, 3, 4, 4, 5, 6, 6, 6, 8, 9, 10, 12, 12, 12, 14, 20,
          24, 24, 24, 30, 48, 50, 60, 60, 72, 120, 192, 600, 1920, 604800 ] )
\end{verbatim}

We see that $S$ has $8\,600$ regular orbits on (the dual space of) $N$.

\section{$G/N \cong J_2$ and $|N| = 2^{28}$}

The group $S = J_2$ has exactly one irreducible $28$-dimensional module
over the field with two elements, up to isomorphism.
This module can be obtained from any of the two absolutely irreducible
$14$-dimensional $S$-modules in characteristic two,
by regarding it as a module over the prime field $\FF_2$.

\begin{verbatim}
    gap> p:= 2;;  d:= 28;;
    gap> t:= CharacterTable( "J2" ) mod p;
    BrauerTable( "J2", 2 )
    gap> irr:= Filtered( Irr( t ), x -> x[1] <= d );;
    gap> Display( t, rec( chars:= irr, powermap:= false, centralizers:= false ) );
    J2mod2
    
           1a 3a 3b 5a 5b  5c  5d 7a 15a 15b
    
    Y.1     1  1  1  1  1   1   1  1   1   1
    Y.2     6 -3  .  A *A   C  *C -1   D  *D
    Y.3     6 -3  . *A  A  *C   C -1  *D   D
    Y.4    14  5 -1  B *B  -C -*C  .   .   .
    Y.5    14  5 -1 *B  B -*C  -C  .   .   .
    
    A = -2*E(5)-2*E(5)^4
      = 1-Sqrt(5) = 1-r5
    B = -3*E(5)-3*E(5)^4
      = (3-3*Sqrt(5))/2 = -3b5
    C = E(5)+2*E(5)^2+2*E(5)^3+E(5)^4
      = (-3-Sqrt(5))/2 = -2-b5
    D = E(5)+E(5)^4
      = (-1+Sqrt(5))/2 = b5
    gap> List( irr, x -> SizeOfFieldOfDefinition( x, p ) );
    [ 2, 4, 4, 4, 4 ]
\end{verbatim}

We use the same approach as in the previous example.

\begin{verbatim}
    gap> tom:= TableOfMarks( "J2" );;
    gap> info:= OneAtlasGeneratingSetInfo( "J2", Dimension, 14, Ring, GF(4) );;
    gap> gens:= List( AtlasGenerators( info ).generators,
    >                 x -> BlownUpMat( Basis(GF(4)), x ) );;
    gap> orbits_from_tom( tom, gens, p );
    rec( decomp := [ 235, 33, 282, 38, 0, 0, 6, 31, 36, 0, 0, 0, 3, 66, 9, 0, 0,
          0, 0, 0, 0, 2, 18, 0, 0, 1, 0, 0, 15, 0, 0, 0, 6, 0, 0, 0, 0, 0, 0, 0,
          12, 0, 0, 5, 0, 1, 0, 0, 0, 3, 0, 0, 0, 0, 0, 0, 0, 0, 0, 0, 3, 1, 3,
          0, 9, 0, 0, 0, 0, 0, 0, 6, 0, 0, 0, 0, 0, 0, 0, 0, 0, 3, 0, 0, 0, 0, 0,
          0, 0, 0, 0, 0, 1, 0, 0, 0, 0, 0, 3, 0, 0, 0, 0, 0, 3, 0, 0, 0, 6, 0, 0,
          0, 0, 0, 0, 9, 0, 0, 0, 0, 0, 0, 0, 0, 0, 1, 0, 0, 0, 0, 1, 0, 0, 0, 0,
          0, 0, 0, 3, 0, 0, 0, 3, 0, 0, 1 ],
      fixed := [ 268435456, 65536, 65536, 65536, 256, 1024, 4096, 1024, 1024,
          256, 256, 256, 64, 1024, 64, 256, 16, 16, 64, 64, 64, 256, 256, 64, 16,
          16, 64, 64, 64, 64, 16, 16, 1024, 4, 4, 4, 4, 16, 16, 16, 64, 16, 16,
          16, 16, 64, 16, 16, 16, 64, 16, 16, 16, 16, 4, 16, 16, 16, 16, 1, 64,
          4, 16, 4, 64, 4, 16, 4, 16, 1, 4, 16, 4, 4, 16, 16, 4, 4, 16, 1, 4, 16,
          1, 1, 1, 16, 4, 16, 16, 16, 1, 4, 16, 4, 4, 1, 4, 64, 4, 4, 4, 16, 4,
          4, 4, 1, 1, 4, 16, 1, 4, 1, 4, 1, 4, 16, 4, 4, 1, 1, 1, 1, 4, 1, 1, 16,
          1, 1, 4, 1, 4, 4, 1, 4, 1, 1, 4, 1, 4, 1, 1, 1, 4, 1, 1, 1 ],
      nonzeropos := [ 1, 2, 3, 4, 7, 8, 9, 13, 14, 15, 22, 23, 26, 29, 33, 41,
          44, 46, 50, 61, 62, 63, 65, 72, 82, 93, 99, 105, 109, 116, 126, 131,
          139, 143, 146 ],
      staborders := [ 1, 2, 2, 3, 4, 4, 4, 6, 6, 6, 8, 8, 9, 10, 12, 12, 14, 16,
          16, 24, 24, 24, 24, 30, 40, 50, 60, 72, 96, 120, 192, 240, 600, 1920,
          604800 ] )
\end{verbatim}

We see that $S$ has $235$ regular orbits on (the dual space of) $N$.

\section{$G/N \cong {}^3D_4(2)$ and $|N| = 2^{26}$}

The group $S = {}^3D_4(2)$ has exactly one irreducible $26$-dimensional module
over the field with two elements, up to isomorphism.
This module is in fact absolutely irreducible.

\begin{verbatim}
    gap> p:= 2;;  d:= 26;;
    gap> t:= CharacterTable( "3D4(2)" ) mod p;
    BrauerTable( "3D4(2)", 2 )
    gap> irr:= Filtered( Irr( t ), x -> x[1] <= d );;
    gap> Display( t, rec( chars:= irr, powermap:= false, centralizers:= false ) );
    3D4(2)mod2
    
           1a 3a 3b 7a 7b 7c 7d 9a 9b 9c 13a 13b 13c 21a 21b 21c
    
    Y.1     1  1  1  1  1  1  1  1  1  1   1   1   1   1   1   1
    Y.2     8  2 -1  A  C  B  1  D  F  E   G   I   H   J   L   K
    Y.3     8  2 -1  B  A  C  1  E  D  F   H   G   I   K   J   L
    Y.4     8  2 -1  C  B  A  1  F  E  D   I   H   G   L   K   J
    Y.5    26 -1 -1  5  5  5 -2  2  2  2   .   .   .  -1  -1  -1
    
    A = 3*E(7)^2+E(7)^3+E(7)^4+3*E(7)^5
    B = 3*E(7)+E(7)^2+E(7)^5+3*E(7)^6
    C = E(7)+3*E(7)^3+3*E(7)^4+E(7)^6
    D = -E(9)^2+E(9)^3-2*E(9)^4-2*E(9)^5+E(9)^6-E(9)^7
    E = -E(9)^2+E(9)^3+E(9)^4+E(9)^5+E(9)^6-E(9)^7
    F = 2*E(9)^2+E(9)^3+E(9)^4+E(9)^5+E(9)^6+2*E(9)^7
    G = E(13)+E(13)^2+E(13)^3+E(13)^5+E(13)^8+E(13)^10+E(13)^11+E(13)^12
    H = E(13)+E(13)^4+E(13)^5+E(13)^6+E(13)^7+E(13)^8+E(13)^9+E(13)^12
    I = E(13)^2+E(13)^3+E(13)^4+E(13)^6+E(13)^7+E(13)^9+E(13)^10+E(13)^11
    J = E(7)^3+E(7)^4
    K = E(7)^2+E(7)^5
    L = E(7)+E(7)^6
\end{verbatim}

We try the same approach as in the examples about the group $J_2$.

\begin{verbatim}
    gap> tom:= TableOfMarks( "3D4(2)" );
    TableOfMarks( "3D4(2)" )
    gap> StandardGeneratorsInfo( tom );
    [ rec( ATLAS := true, description := "|z|=8, z^4=a, |b|=9, |ab|=13, |abb|=8", 
          generators := "a, b", 
          script := [ [ 1, 8, 4 ], [ 2, 9 ], [ 1, 1, 2, 1, 13 ], 
              [ 1, 1, 2, 1, 2, 1, 8 ] ], standardization := 1 ) ]
    gap> info:= OneAtlasGeneratingSetInfo( "3D4(2)", Dimension, 26, Ring, GF(2) );;
    gap> gens:= AtlasGenerators( info ).generators;;
    gap> map:= GroupGeneralMappingByImages( UnderlyingGroup( tom ),
    >      Group( gens ), GeneratorsOfGroup( UnderlyingGroup( tom ) ), gens );;
    gap> IsGroupHomomorphism( map );
    true
\end{verbatim}

Now we apply the function \verb|orbits_from_tom|.

\begin{verbatim}
    gap> orbsinfo:= orbits_from_tom( tom, gens, p );;
    gap> orbsinfo.fixed[1];
    67108864
    gap> orbsinfo.decomp[1];
    0
\end{verbatim}

Unfortunately, $S$ has no regular orbit on (the dual of) $N$.
However, there is one orbit whose point stabilizer in $S$ is a dihedral group
$D_{18}$ of order $18$.

\begin{verbatim}
    gap> orbsinfo.staborders;
    [ 16, 16, 18, 42, 48, 52, 64, 72, 392, 1008, 1536, 3024, 3072, 3584, 258048, 
      211341312 ]
    gap> orbsinfo.nonzeropos[3];
    446
    gap> orbsinfo.decomp[446];
    1
    gap> u:= RepresentativeTom( tom, 446 );
    <permutation group of size 18 with 2 generators>
    gap> IsDihedralGroup( u );
    true
\end{verbatim}

Thus there ia a linear character $\lambda$ of $N$ whose inertia subgroup
$T = I_G(\lambda)$ has the structure $N.D_{18}$.
Now $\Irr( T | \lambda )$ can be identified with those irreducibles of
$T/\ker(\lambda)$ that restrict nontrivially to $N/\ker(\lambda)$,
and there are only two groups, up to isomorphism, that can occur as
$T/\ker(\lambda)$.

\begin{verbatim}
    gap> cand:= Filtered( AllSmallGroups( 36 ),
    >             x -> Size( Centre( x ) ) = 2 and
    >                  IsDihedralGroup( x / Centre( x ) ) );
    [ <pc group of size 36 with 4 generators>, 
      <pc group of size 36 with 4 generators> ]
    gap> List( cand, StructureDescription );
    [ "C9 : C4", "D36" ]
\end{verbatim}

These two groups are a split and a nonsplit extension of the cyclic group
of order $18$ with a group of order two that acts by inverting.
In other words, these two groups are the direct product of $D_{18}$ with
a cyclic group of order two and the subdirect product of $D_{18}$ with
a cyclic group of order four.

Both groups possess irreducible characters of degree two, one rational
valued and the other not, which restrict nontrivially to the centre.

\begin{verbatim}
    gap> Display( CharacterTable( "Dihedral", 18 ) );
    Dihedral(18)
    
         2  1  .  .  .  .  1
         3  2  2  2  2  2  .
    
           1a 9a 9b 3a 9c 2a
        2P 1a 9b 9c 3a 9a 1a
        3P 1a 3a 3a 1a 3a 2a
    
    X.1     1  1  1  1  1  1
    X.2     1  1  1  1  1 -1
    X.3     2  A  B -1  C  .
    X.4     2  B  C -1  A  .
    X.5     2 -1 -1  2 -1  .
    X.6     2  C  A -1  B  .
    
    A = -E(9)^2-E(9)^4-E(9)^5-E(9)^7
    B = E(9)^2+E(9)^7
    C = E(9)^4+E(9)^5
\end{verbatim}

By \cite[Lemma~5.1~(ii)]{DNT}, we are done.

\section{$G/N \cong {}^3D_4(2)$ and $|N| = 3^{25}$}

The group $S = {}^3D_4(2)$ has exactly one irreducible $25$-dimensional module
over the field with three elements, up to isomorphism.
This module is in fact absolutely irreducible.

\begin{verbatim}
    gap> p:= 3;;  d:= 25;;
    gap> t:= CharacterTable( "3D4(2)" ) mod p;
    BrauerTable( "3D4(2)", 3 )
    gap> irr:= Filtered( Irr( t ), x -> x[1] <= d );;
    gap> Display( t, rec( chars:= irr, powermap:= false, centralizers:= false ) );
    3D4(2)mod3
    
           1a 2a 2b 4a 4b 4c 7a 7b 7c 7d 8a 8b 13a 13b 13c 14a 14b 14c 28a 28b 28c
    
    Y.1     1  1  1  1  1  1  1  1  1  1  1  1   1   1   1   1   1   1   1   1   1
    Y.2    25 -7  1  5 -3  1  4  4  4 -3 -1 -1  -1  -1  -1   .   .   .  -2  -2  -2
\end{verbatim}

We use the same approach as in the examples about the group $J_2$.

\begin{verbatim}
    gap> tom:= TableOfMarks( "3D4(2)" );;
    gap> info:= OneAtlasGeneratingSetInfo( "3D4(2)", Dimension, d, Ring, GF(p) );;
    gap> gens:= AtlasGenerators( info ).generators;;
    gap> orbsinfo:= orbits_from_tom( tom, gens, p );;
    gap> orbsinfo.fixed[1];
    847288609443
    gap> orbsinfo.decomp[1];
    3551
\end{verbatim}

We see that $S$ has $3\,551$ regular orbits on (the dual space of) $N$.

\bibliographystyle{amsalpha}
\newcommand{\etalchar}[1]{$^{#1}$}
\providecommand{\bysame}{\leavevmode\hbox to3em{\hrulefill}\thinspace}
\providecommand{\MR}{\relax\ifhmode\unskip\space\fi MR }
\providecommand{\MRhref}[2]{%
  \href{http://www.ams.org/mathscinet-getitem?mr=#1}{#2}
}
\providecommand{\href}[2]{#2}


\end{document}